\documentclass[12pt]{article}
\usepackage{graphicx}

\usepackage{amsmath,amssymb}
\usepackage{mathrsfs}
\usepackage{amsfonts}
\usepackage{latexsym}
\usepackage{amsthm}
\usepackage[T2A]{fontenc}     
\usepackage[cp1251]{inputenc} 
\usepackage{graphicx}

\newcommand{\er}[1]{{\rm(\ref{#1})}}
\def\lb{\label}
\theoremstyle{plain}
\newtheorem{theorem}{\bf Theorem}[section]
\newtheorem{lemma}[theorem]{\bf Lemma}
\theoremstyle{remark}

\newtheorem{proposition}[theorem]{\bf Proposition}

\begin{document}

\def\a{\alpha} \def\cA{{\cal A}} \def\bA{{\bf A}}  \def\mA{{\mathscr A}}
\def\b{\beta}  \def\cB{{\cal B}} \def\bB{{\bf B}}  \def\mB{{\mathscr B}}
\def\g{\gamma} \def\cC{{\cal C}} \def\bC{{\bf C}}  \def\mC{{\mathscr C}}
\def\G{\Gamma} \def\cD{{\cal D}} \def\bD{{\bf D}}  \def\mD{{\mathscr D}}
\def\d{\delta} \def\cE{{\cal E}} \def\bE{{\bf E}}  \def\mE{{\mathscr E}}
\def\D{\Delta} \def\cF{{\cal F}} \def\bF{{\bf F}}  \def\mF{{\mathscr F}}
\def\c{\chi}   \def\cG{{\cal G}} \def\bG{{\bf G}}  \def\mG{{\mathscr G}}
\def\z{\zeta}  \def\cH{{\cal H}} \def\bH{{\bf H}}  \def\mH{{\mathscr H}}
\def\e{\eta}   \def\cI{{\cal I}} \def\bI{{\bf I}}  \def\mI{{\mathscr I}}
\def\p{\psi}   \def\cJ{{\cal J}} \def\bJ{{\bf J}}  \def\mJ{{\mathscr J}}
\def\vT{\Theta}\def\cK{{\cal K}} \def\bK{{\bf K}}  \def\mK{{\mathscr K}}
\def\k{\kappa} \def\cL{{\cal L}} \def\bL{{\bf L}}  \def\mL{{\mathscr L}}
\def\l{\lambda}\def\cM{{\cal M}} \def\bM{{\bf M}}  \def\mM{{\mathscr M}}
\def\L{\Lambda}\def\cN{{\cal N}} \def\bN{{\bf N}}  \def\mN{{\mathscr N}}
\def\m{\mu}    \def\cO{{\cal O}} \def\bO{{\bf O}}  \def\mO{{\mathscr O}}
\def\n{\nu}    \def\cP{{\cal P}} \def\bP{{\bf P}}  \def\mP{{\mathscr P}}
\def\r{\rho}   \def\cQ{{\cal Q}} \def\bQ{{\bf Q}}  \def\mQ{{\mathscr Q}}
\def\s{\sigma} \def\cR{{\cal R}} \def\bR{{\bf R}}  \def\mR{{\mathscr R}}
\def\S{\Sigma} \def\cS{{\cal S}} \def\bS{{\bf S}}  \def\mS{{\mathscr S}}
\def\t{\tau}   \def\cT{{\cal T}} \def\bT{{\bf T}}  \def\mT{{\mathscr T}}
\def\f{\phi}   \def\cU{{\cal U}} \def\bU{{\bf U}}  \def\mU{{\mathscr U}}
\def\F{\Phi}   \def\cV{{\cal V}} \def\bV{{\bf V}}  \def\mV{{\mathscr V}}
\def\P{\Psi}   \def\cW{{\cal W}} \def\bW{{\bf W}}  \def\mW{{\mathscr W}}
\def\o{\omega} \def\cX{{\cal X}} \def\bX{{\bf X}}  \def\mX{{\mathscr X}}
\def\x{\xi}    \def\cY{{\cal Y}} \def\bY{{\bf Y}}  \def\mY{{\mathscr Y}}
\def\X{\Xi}    \def\cZ{{\cal Z}} \def\bZ{{\bf Z}}  \def\mZ{{\mathscr Z}}
\def\O{\Omega}
\def\ve{\varepsilon}
\def\vt{\vartheta}
\def\vp{\varphi}
\def\vk{\varkappa}

\def\mM{M}
\def\mB{B}
\def\mR{R}

\def\mA{{\mathscr A}}
\def\mB{{\mathscr B}}
\def\mC{{\mathscr C}}
\def\mD{{\mathscr D}}
\def\mE{{\mathscr E}}
\def\mF{{\mathscr F}}
\def\mG{{\mathscr G}}
\def\mH{{\mathscr H}}
\def\mI{{\mathscr I}}
\def\mJ{{\mathscr J}}
\def\mK{{\mathscr K}}
\def\mL{{\mathscr L}}
\def\mM{{\mathscr M}}
\def\mN{{\mathscr N}}
\def\mO{{\mathscr O}}
\def\mP{{\mathscr P}}
\def\mQ{{\mathscr Q}}
\def\mR{{\mathscr R}}
\def\mS{{\mathscr S}}
\def\mT{{\mathscr T}}
\def\mU{{\mathscr U}}
\def\mV{{\mathscr V}}
\def\mW{{\mathscr W}}
\def\mX{{\mathscr X}}
\def\mY{{\mathscr Y}}
\def\mZ{{\mathscr Z}}

\newcommand{\gA}{\mathfrak{A}}
\newcommand{\gB}{\mathfrak{B}}
\newcommand{\gC}{\mathfrak{C}}
\newcommand{\gD}{\mathfrak{D}}
\newcommand{\gE}{\mathfrak{E}}
\newcommand{\gF}{\mathfrak{F}}
\newcommand{\gG}{\mathfrak{G}}
\newcommand{\gH}{\mathfrak{H}}
\newcommand{\gI}{\mathfrak{I}}
\newcommand{\gJ}{\mathfrak{J}}
\newcommand{\gK}{\mathfrak{K}}
\newcommand{\gL}{\mathfrak{L}}
\newcommand{\gM}{\mathfrak{M}}
\newcommand{\gN}{\mathfrak{N}}
\newcommand{\gO}{\mathfrak{O}}
\newcommand{\gP}{\mathfrak{P}}
\newcommand{\gR}{\mathfrak{R}}
\newcommand{\gS}{\mathfrak{S}}
\newcommand{\gT}{\mathfrak{T}}
\newcommand{\gU}{\mathfrak{U}}
\newcommand{\gV}{\mathfrak{V}}
\newcommand{\gW}{\mathfrak{W}}
\newcommand{\gX}{\mathfrak{X}}
\newcommand{\gY}{\mathfrak{Y}}
\newcommand{\gZ}{\mathfrak{Z}}

\def\Z{{\Bbb Z}}
\def\R{{\Bbb R}}
\def\C{{\Bbb C}}
\def\T{{\Bbb T}}
\def\N{{\Bbb N}}
\def\S{{\Bbb S}}
\def\H{{\Bbb H}}
\def\J{{\Bbb J}}
\def\dD{{\Bbb D}}

\def\qqq{\qquad}
\def\qq{\quad}
\newcommand{\ma}{\begin{pmatrix}}
\newcommand{\am}{\end{pmatrix}}
\newcommand{\ca}{\begin{cases}}
\newcommand{\ac}{\end{cases}}
\let\ge\geqslant
\let\le\leqslant
\let\geq\geqslant
\let\leq\leqslant
\def\ma{\left(\begin{array}{cc}}
\def\am{\end{array}\right)}
\def\iint{\int\!\!\!\int}
\def\lt{\biggl}
\def\rt{\biggr}
\let\geq\geqslant
\let\leq\leqslant
\def\[{\begin{equation}}
\def\]{\end{equation}}
\def\wh{\widehat}
\def\wt{\widetilde}
\def\pa{\partial}
\def\sm{\setminus}
\def\es{\emptyset}
\def\no{\noindent}
\def\ol{\overline}
\def\iy{\infty}
\def\ev{\equiv}
\def\/{\over}
\def\ts{\times}
\def\os{\oplus}
\def\ss{\subset}
\def\h{\hat}
\def\Ra{\Rightarrow}
\def\Re{\mathop{\rm Re}\nolimits}
\def\Im{\mathop{\rm Im}\nolimits}
\def\supp{\mathop{\rm supp}\nolimits}
\def\sign{\mathop{\rm sign}\nolimits}
\def\Ran{\mathop{\rm Ran}\nolimits}
\def\Ker{\mathop{\rm Ker}\nolimits}
\def\Tr{\mathop{\rm Tr}\nolimits}
\def\const{\mathop{\rm const}\nolimits}
\def\dist{\mathop{\rm dist}\nolimits}
\def\diag{\mathop{\rm diag}\nolimits}
\def\Wr{\mathop{\rm Wr}\nolimits}
\def\Iso{\mathop{\rm Iso}\nolimits}
\def\Res{\mathop{\rm Res}\nolimits}
\def\BBox{\hspace{1mm}\vrule height6pt width5.5pt depth0pt \hspace{6pt}}
\def\Diag{\mathop{\rm Diag}\nolimits}

\def\Twelve{
\font\Tenmsa=msam10 scaled 1200 \font\Sevenmsa=msam7 scaled 1200
\font\Fivemsa=msam5 scaled 1200 \textfont\msbfam=\Tenmsb
\scriptfont\msbfam=\Sevenmsb \scriptscriptfont\msbfam=\Fivemsb

\font\Teneufm=eufm10 scaled 1200 \font\Seveneufm=eufm7 scaled 1200
\font\Fiveeufm=eufm5 scaled 1200
\textfont\eufmfam=\Teneufm \scriptfont\eufmfam=\Seveneufm
\scriptscriptfont\eufmfam=\Fiveeufm}

\def\Ten{
\textfont\msafam=\tenmsa \scriptfont\msafam=\sevenmsa
\scriptscriptfont\msafam=\fivemsa

\textfont\msbfam=\tenmsb \scriptfont\msbfam=\sevenmsb
\scriptscriptfont\msbfam=\fivemsb

\textfont\eufmfam=\teneufm \scriptfont\eufmfam=\seveneufm
\scriptscriptfont\eufmfam=\fiveeufm}

\title {Inverse resonance scattering for Jacobi operators}

\author{
 Evgeny Korotyaev
\begin{footnote}
{School of Math., Cardiff Univ., Senghennydd Road, Cardiff, CF24
4AG, UK, e-mail: KorotyaevE@cf.ac.uk}
\end{footnote}
}

\maketitle

\begin{abstract}
\no    We consider the Jacobi operator $(Jf)_n=
a_{n-1}f_{n-1}+a_nf_{n+1}+b_nf_n$ on $\Z$ with a real compactly
supported sequences $(a_n-1)_{n\in\Z}$ and $(b_n)_{n\in\Z}$. We give
the solution of two inverse problems (including characterization): $
(a,b)\to \{$zeros of the reflection coefficient$\}$ and $(a,b)\to
\{$bound states and resonances$\}$. We describe the set of
"iso-resonance operators $J$", i.e., all operators $J$ with the same
resonances and bound states.
\end{abstract}




\section {Introduction}
\setcounter{equation}{0}

We consider the Jacobi operator $J$ acting on $\ell^2$ and given by
$$
 (Jf)_n= a_{n-1}f_{n-1}+a_nf_{n+1}+b_nf_n,\ \ n\in\Z,\qqq f=(f_n)_{-\iy}^{\iy}
 \in \ell^2(\Z),
$$
where a real compactly supported sequence  $q=(q_n)_{n\in \Z},\
q_{2n-1}=b_n, \ q_{2n}=1-a_n$ satisfies
\begin{multline}
 \lb{qs}
q=(q_n)_{n\in \Z}\in \gX_\n^\t=\gX_\n^\t(p)=\ell_{1+\n,2p-\t},\qq
for\  some \ \n,\t\in \{0,1\},\qq p\in \N,
\\
\ell_{m,k}=\rt\{(q_n)_{n\in \Z}\in \ell^2(\Z): q_m\ne 0,
q_k\ne 0,\ q_{2s}<1, \ q_n=0, \ all\ s\in \Z,\ n\in \Z\sm [m,k]\rt\},
\end{multline}
for some $m,k\in \Z$. For fixed $p\ge1$ the sequence $q\in \gX_{0}^0$ has the max support and $q\in \gX_{1}^1$ has the min
support and $q\ne 0$.  It is well known that the spectrum of $J$ has the form
\[
\lb{sJ} \s(J)=\s_{ac}(J)\cup \s_d(J),\qqq \s_{ac}(J)=[-2,2],\qqq
\s_d(J)\ss \R\sm [-2,2].
\]
Define the new variable $z\in \dD_1=\{z\in
\C:|z|<1\}$ by $\l=\l(z)=z+{1\/z}$. Here $\l(z)$ is a conformal mapping from $\dD_1$ onto
$\C\sm [-2,2]$. Denote by $\p^\pm =(\p_n^\pm(z))_{-\iy}^\iy$ the
fundamental solutions to
\[
\lb{je}
a_{n-1}\p^\pm _{n-1}+a_n\p^\pm _{n+1}+b_n\p^\pm _n=(z+{1\/z})\p^\pm _n,\ \
n\in\Z,
\]
\[
\lb{4}
\p_n^+(z)=z^{n},\ n>p,\qqq  and \qqq \p_{n}^-(z)=z^{-n},\
n\le 0, \qqq |z|\le 1.
\]
Define the Wronskian $\{f,u\}_n=a_n(f_nu_{n+1}-u_nf_{n+1})$
for sequences $u=(u_n)_{-\iy}^{\iy}, f=(f_n)_{-\iy}^{\iy}$.
If $f,u$ are some solutions of \er{je}, then
$\{f,u\}_n$ does not depend on $n$.
The following identities hold true:
\[
\lb{5} \p^+= A\wt\p^-+B\p^-, \qqq on \qq  \S^1_0=\S^1\sm \{\pm 1\},\ \  \qq
\wt\p^\pm=\p^\pm(z^{-1}),
\]
where $\S^1=\{z\in\C: |z|=1\},$ and $\wt v=v(z^{-1})$ for a function $v(z)$ and
\[
\lb{6} A={w\/1-z^2},\qq  w=z\{\p^+, \p^-\}_n,\qqq B={z^2s\/\e},\qq
s={\{\p^+, \wt\p^-\}_n\/z^2},\qqq \e=z-{1\/z}.
\]
Note that if $q=0$, then $w=1-z^2, A=1, B=0$.  The S-matrix is given by
\[
\lb{sm1} \cS_M(z)=\ma A(z)^{-1}& R_-(z)\\ R_+(z)&A(z)^{-1}\am ,  \qq z\in \S^1,
\qqq R_-={B\/A}=-{z^3s\/w}, \qq R_+=-{\wt B\/A}={\wt s\/zw},
\]
where ${1\/A}$ is the transmission coefficient and $R_{\pm}$ is the
reflection coefficient. For each $z\in \S^1_0$ the scattering matrix $\cS_M$ is unitary and
 satisfies
\[
\lb{7} |A(z)|^2 =1+|B(z)|^2\qq \Leftrightarrow \qq  w(z)w(z^{-1})+\e^2(z)=s(z)s(z^{-1}),
\]
\[
\lb{sets} \det \cS_M(z) ={\ol A(z)\/A(z)}=-z^2{\ol
w(z)\/w(z)}=-z^2{w(z^{-1})\/w(z)}.
\]
In the case of compactly supported $q$
it is convenient to work with the  polynomials $w,s$.
Thus $\cS_M$ is meromorphic in the complex plane $\C$ and the
poles of $\cS_M$ are the zeros of $w$. Let $\l_j,j\in
\Z_N=\{n_-, ..,-1,1,...,n_+\}, N=n_+-n_-\ge 0$ be the bound states of
$J$ and $\l_j=z_j+z_j^{-1}$, where $z_j$ are all zeros of $w$ in
$\dD_1$ and the sequence $E_N=(z_j)_{j\in\Z_N}$ satisfies
\[
\lb{eg1} -1<z_{n_-}<...<z_{-1}<0<z_{1}<...<z_{n_+}<1\qq for \qq
some \ \pm n_\pm\ge 0.
\]
Recall that $w$ has only these real zeros  in $\dD_1$ (otherwise $J$ has a
non-real eigenvalue).

The main goal of this paper is to prove the following results:

\no {\it i) the mapping $(a,b)\to \{$bound states,  zeros of the reflection
coefficient$\}$ is 1-to-1 and onto,

\no ii)  the mapping $(a,b)\to \{$bound states, resonances and
some sequence $\s=(\s_n)_0^m, \s_n=\pm 1,  \}$ is 1-to-1 and
onto,

\no iii) we characterize "iso-resonance Jacobi operators $J$", i.e., all operators $J$ with the same resonances and bound states.}

The inverse spectral problem consists of the following parts:\\
i) Uniqueness.
Prove that the spectral data uniquely determine the potential.\\
ii) Characterization.
Give conditions for some data to be the spectral data of some potential.\\
iii) Reconstruction. Give an  algorithm for recovering the
potential from the spectral data.\\
 iv) A priori estimates. Obtain estimates of the potential in
 terms of the  spectral data.

 We define the class of scattering data $(s, E_N)$, where $s(z)$ is a polynomial associated with the zeros of the reflection coefficient  and $E_N=(z_j)_{j\in\Z_N}$ is a sequence  of  all zeros  of $w$ in $\dD_1$.

 \no {\bf Definition S.} \no {\it By
$\cS_\n^\t=\cS_\n^\t(p), \n,\t\in \{0,1\}, p\in \N$ we will denote
the class of $(f,K_f^N)$, for some $N\ge 0$, where $f$ is a real polynomial
  (i.e., a polynomial with real coefficients) given by
$$
f=Cz^\n\prod_{1}^{m}(z-\z_n),\qq  where \ C\in\R\sm\{0\},\ \ \z_n
\in \C\sm\{0\},\qq         m=2p-1-\t-\n,
$$
$K_f^N=(z_j)_{j\in\Z_N}$ is a sequence  of zeros of the
function $f(z)f(z^{-1})-\e^2(z)$, which  satisfy
\[
\lb{S1} -1<z_{n_-}<...<z_{-1}<0<z_{1}<...<z_{n_+}<1,\qq \pm n_\pm\ge 0,
\qq N=n_+-n_-,
\]
\[
\lb{S2} (-1)^{j}z_jf(z_j^{\pm 1})>0,\ \ j\in \Z_N,\ and \qqq
(-1)^{n_\pm}f(\pm 1)\ge 0\ .
\]
}
\no  {\bf Remark}. i) Below we will show that if $q\in \gX_\n^\t$,
then  $(s,E_N)\in \cS_\n^\t$, where $E_N=K_s^N$ is a
sequence of all zeros of $w$ in $\dD_1$. Note that $K_s^0=\es$. ii)
It is possible that the function $s(z)s(z^{-1})-\e^2(z)$ has $N_1>N$
zeros on $(-1,1)$, and it is important that we take some of
them with the needed properties \er{S2}. In particular, we emphasize
that $E_N$ is not uniquely determined by $s$.

  \begin{theorem} \lb{T1}
  Let $\t,\n\in \{0,1\}$ and $m\ge 3$.
The mapping $\gJ: \gX_\n^\t\to \cS_\n^\t$ given by $\gJ (q)=(s,E_N)$
is  one-to-one and onto.
\end{theorem}

\no {\bf Remark.} 1) In the proof of Theorem \ref{T1} we present
algorithm for recovering the potential $q$ from the spectral data
$(s,E_N)$. The potential $q$ is uniquely determined by the Marchenko
equations \er{GLM1}-\er{recab} (in terms of $(s,E_N)$). Here we use
the results about the Marchenko equation  from \cite{Te}.
Standard spectral data for the inverse problem for the Jacobi operator $J$ are $(s, z_j, m_j^\pm, j\in \Z_N)$, where
$m_j^\pm$ is so-called norming constant  given by
\[
m_j^\pm =\sum_{n\in \Z} \p_n^\pm (z_j)^2,\qqq j\in\Z_N,
\]
see \cite{Te}. In Theorem \ref{T1} instead of the norming constants
we need the condition \er{S2}.

\no 2) In Section 3 we show simple examples of the scattering data for the case $p=1,2$.

\no 3)  We briefly indicate how to prove  Theorem 1.1. Firstly, we
show that if $q\in \gX_\t^\n$, then $(s,E_N)\in \cS_\n^\t$. Here we
check condition \er{S2}. Secondly, we consider the inverse mapping.
Suppose $(s,K_s^N)\in \cS_\t^\n$. Then in order to determine $w$ we solve
the functional equation \er{7} in some class of
polynomials (see  Theorem \ref{Ta1})  and this gives the reflection
coefficient $R_\pm$. Thirdly, \er{sw3-2}, \er{sw3-3}
yield the norming constants $m_j^\pm$. Then we check that $R_\pm, z_j, m_j^\pm,
j\in\Z_N$ satisfy conditions from Theorem \ref{Tist} of Teschl \cite{Te},
which we recall for the sake of the reader in Sect.3. Then we obtain a bijection
of our mapping.

\no 4) Assume that $s(\cdot)=0$ for some $q\in \gX_\n^\t$. Then
\er{7} yields $w(z)w(z^{-1})=-\t(z)^2, z\in\C\sm\{0\}$. It is
impossible since $w$ is a polynomial and we have
a contradiction. In fact we deduce that if $s(\cdot)=0$ for some
compactly supported "potential" $q$, then $q=0$.

 Let $\# (f, I)$ denote the number of zeros of a function $f$ on the set $I$.

\no {\bf Definition W.} \no {\it By $\cW_\n^\t=\cW_\n^\t(p), \n,\t\in \{0,1\}, p\in \N$ we will denote the class of polynomials
\[
\lb{defW1} w=C_w\prod_1^{m}(z-\r_n), \qq C_w\in \R\sm\{0\}, \qq
m=2p-1-\t-\n ,\qq\r_n \in \C,\qq w(0)>0,
\]

\no i) $w$ is real on $\R$ and  $|w(z)|\ge |\e(z)|$ for  any $|z|=1$, where $\e=z-{1\/z}$,

\no ii)  $w$ has only real simple zeros $z_j, j\in \Z_N$ in $\ol\dD_1$,
and the sequence $E_N=(z_j)_{j\in \Z_N}$ and the function $F(z)=w(z)w(z^{-1})+\e^2(z)$
satisfies}
\begin{multline}
\lb{defW2}
-1<z_{n_-}<...<z_{-1}<0<z_{1}<...<z_{n_+}<1,\qq for \ some\ \pm n_\pm \ge 0, \ N=n_+-n_-\ge 0,
\\
{1\/2}\#(F,(z_{n_\pm},z_{n_\pm}^{-1}))=even \ge 0,\qq
\#(F,I_j)=even \ge 2,\qq
\\
I_0=(z_{-1},z_{1}), \qq  I_j=(z_j,z_{j+1}),\qq
j\in \Z_N\sm \{-1,n_{n_+}\}=\{n_-,...,-2,1,2,..,n_+-1\}.
\end{multline}

 We describe the properties of $s,w$
and the sequence of zeros  $E_N=(z_j)_{j\in\Z_N}\ss (-1,1)^N$.

\begin{proposition}
\lb{sw3}
{\it Let $q\in \gX_\n^\t$ for some $\t,\n\in \{0,1\}$.  Then  $(s,E_N)\in \cS_\n^\t, w\in \cW_\n^\t$ and for each $j\in \Z_N$ the following identities hold true
\[
\lb{sw3-1} \p^+(z_j)=B(z_j)\p^-(z_j),\qqq  B(z_j^{-1})B(z_j)=-1,\qqq
s(z_j^{-1})s(z_j)=\e^2(z_j),
\]
\[
\lb{sw3-2}
 m_j^+=-z_jA'(z_j)B(z_j)={z_j^2\/\e^2(z_j)}w'(z_j)s(z_j),\qqq
\]
\[
\lb{sw3-3} m_j^+=m_j^-B^2(z_j),
\]
\[
\lb{sw3-4} z_j(-1)^{j}s(z_j^{\pm 1})>0,\qqq (-1)^{j}z_jw'(z_j)>0,
\]
\[
\lb{sw3-5} (-1)^{n_\pm}w(\pm 1)\ge 0,\qqq  s(\pm 1)=\pm w(\pm 1).
\]
Moreover, let
 $w(z)=\sum_1^{2p}\check w_n z^{n-1},\  \ s(z)=\sum_1^{2p}\check s_nz^{n-1},\qq \check  w=(\check  w_n)_1^{2p},
 \check  s=(\check  s_n)_1^{2p}\in \R^{2p}$ and let $Vh=(0,h_1,..,h_{2p-1}), h=(h_n)_1^{2p}$.
Then}
\[
\lb{trp}
2+(\check  s,\check  s)=(\check  w,\check  w),\qq  (V^2\check  s,\check s)=1+(V^2\check
w,\check  w)   \qq (V^{2k-1}\check  s,\check s)=(V^{2k-1}\check w,\check  w),
k=1,....
\]

\end{proposition}

Below we will sometimes write $w(z,q), s(z,q),..$, instead of $w(z), s(z), ..$, when several potentials are being dealt with. For $q\in \gX_\n^\t$  the iso-resonance set of potentials is given  by
\[
\label{IsoDef} \Iso(q)=\left\{r\in \gX_\n^\t:
w(\cdot,q)=w(\cdot,r)\right\}.
\]
We will describe $ \Iso(q)$.
Assume that we know $w$ and we need to recover the polynomial
$s$. Due to Theorem \ref{T1} the function
$F(z)=w(z)w(z^{-1})+\e^2(z)=s(z)s(z^{-1})$ has the zeros $t_n, t_{n+m}=t_n^{-1},
n=1,..,m$ counted with multiplicity and given by
\begin{multline}
 \lb{zef}
0<|t_1|\le ...\le |t_{m}|\le 1,\qq (t_n)_{1}^{m} \ss U=\dD_1\cup
\ol\S_+^1,\qq \S_+^1=\C_+\cap \S^1,
\\ \arg t_n\in [0,2\pi),\qq and \qq
if\ \ |t_n|=|t_k|,\ \arg t_n\le \arg t_k \Rightarrow n\le k.
\end{multline}
Note that if $|t_n|=1$, then $\Im t_n\ge 0$.
Hence $(t_n)_1^m$ is a uniquely defined sequence of all zeros $\ne
0$ of $F$ in the set $U$. Thus $t_n, t_n^{-1}, n=1,..,m$ are all
zeros of $F$ and $t_n$ or $t_n^{-1}$ is a zero of $s$.

If $q\in \gX_\t^\n$, then $s=C_sz^\n\prod_1^{m}(z-\z_n)$, where
each $\z_n\ne 0, n=1,..,m$ and recall
that $\n\in \{0,1\}$. The sequence $\s=(\s_n)_0^{m}$ is defined by
\[
\lb{setX0}
\s=(\s_n)_0^{m}\in \{\pm 1\}^{m+1},\qq \s_0=\sign C_s,\qq and \qq \z_n=t_n^{\s_n}, \
n=1,..,m.
\]
For each $w\in \cW_\t^\n$ we define a set $\X_w$ of all possible sequences
$\s=(\s_n)_0^m$ by
$$
 \lb{setX}
 \X_g=\rt\{\s=(\s_n)_0^m\ss \{-1,1\}^{m+1}: (s,E_N)\in \cS_\n^\t,
 \  where \  s=Cz^\n\prod_1^m (z-t_n^{\s_n}),\s_0=\sign C\rt\}.
$$
In particular we have:

\no I) {\it If $w(1)=w(-1)=0=N$, then $(\s_n)_0^m$ is any sequence from
$\{-1,1\}^{m+1}$, under the condition that $f$ is real.

\no II) If $w(1)\neq 0$ (or $w(-1)\neq 0$), then $(-1)^{n_+}s(1)>0$
(or $(-1)^{n_-}s(1)<0$ ) gives  $\s_0=\sign C_s$.

\no III) If $N\ge 1$, then condition $(-1)^{j}z_js(z_j)>0$ for some
$j\in \Z_N$ gives $\s_0=\sign C_s$.

\no IV) If $N\ge 2$, then the function $s$ has an odd number
$\ge1$ of zeros on each of the intervals $(z_{n_-}, z_{n_-+1}), ..., (z_{-2}, z_{-1})$
$(z_{-1}, z_{1})$ and $(z_{1}, z_{2}),..., (z_{n_+-1}, z_{n_+})$.}

 Our goal is to show that the
spectral data $\X_w$ give the "proper" parametrization of the set
$\Iso(q)$.  Our main Theorem \ref{T2} shows that
$\s\in\X_w$ are {\bf almost free parameters}. Namely, we prove that
if the function $w(z,q)$ is fixed, then each $\s_n$
can be changed in an almost arbitrary way.

\no \begin{theorem} \lb{T2} Let $\t,\n\in \{0,1\}$ and $m\ge 3$.

\no i)    The mapping $\gJ_{R}: \gX_\n^\t\to \{ (w, \X_w), w\in
\cW_\n^\t\}$ given by $q\to (w,(\s_n)_0^m)$  is one-to-one and onto.

\no  ii) Let $q\in  \gX_\n^\t$. Then the mapping $\P: \Iso(q)\to \X_w$, given by   $r\to \s(r)$
(see \er{setX0}) is a bijection between the set of potentials $r\in \Iso(q)$ and the set
of sequences $\s(r)\in  \X_w, w=w(\cdot,q)$.

\end{theorem}

A great number of papers are devoted to the inverse problem
for the Schr\"odinger operator,
(see a book \cite{M} ,  papers \cite{Fa}, \cite{DT}, \cite{Me} and ref. therein).

A lot of papers are devoted to the resonances for the 1D
Schr\"odinger operator, see \cite{F},
\cite{K1}, \cite{K2}, \cite{K3}, \cite{S},\cite{Z}, \cite{Z1}. We recall that Zworski \cite{Z}
obtained the first results about the distribution of resonances for
the Schr\"odinger operator with compactly supported potentials on
the real line.  Korotyaev obtained the characterization (plus
uniqueness and recovering) of $S$-matrix for the Schr\"odinger
operator  with a compactly supported potential on the real line
\cite{K1} and on the half-line \cite{K2}. In \cite{K3} for the
Schr\"odinger operator on the half line the
stability result was given: if $\vk^0=\{\vk^0\}_1^\iy$ is a sequence of zeros
(eigenvalues and resonances) of the Jost function for some real
compactly supported potential $q_0$ and $\vk-\vk^0\in\ell_\ve^2$ for
some $\ve>1$, then $\vk$ is the sequence of zeros of the Jost
function for some unique real compactly supported potential.

There are a lot of papers and books devoted to the scattering for Jacobi operators,
see \cite{C1}, \cite{C2}, \cite{CC}, \cite{CK}, \cite{G1}, \cite{G2},
\cite{NMPZ}, \cite{Te}, \cite{T1}, \cite{T2}.
In the case of Jacobi operators also there are papers about the
inverse resonance problem, see \cite{BNW}, \cite{MW}, \cite{DS1},
\cite{DS2}. In particular, some progress was made by Damanik and
Simon in \cite{DS1}, \cite{DS2}, where they described the $S$-matrix
for Jacobi operators on the half-lattice both for finite-support and
exponentially decay perturbations.

\section {Proof of main Theorems}
\setcounter{equation}{0}

We recall well-known facts  from \cite{Te}.

\begin{lemma}
\lb{f01} Let $q\in \ell_{1,2p}^2$.  Then each function $\p_{p-n}^+(z),
n=1,2,..,p$ is a real polynomial and satisfies
\[
\lb{f01-1} \p_{p}^+(z)={z^p\/a_p},\qq
 \p_{p-n}^+(z)={z^{p+n}\/\e_n}\rt(c_p-{c_p\b_{p-n+1}+b_pa_p^2\/z}+{O(1)\/z^{2}})\rt),
 \qqq c_n=1-a_n^2,
 \qq
\]
\[
\lb{f01-2}
\p_{2}^+(z)={z^{2p-2}\/\e_2}\rt(c_p-{c_p\b_{3}+b_pa_p^2\/z}+{O(1)\/z^{2}}\rt),
\]
\[
\lb{f01-3}
\p_{1}^+(z)={z^{2p-1}\/\e_1}\rt(c_p-{c_p\b_{2}+b_pa_p^2\/z}+{O(1)\/z^{2}}\rt)
\]
as $z\to \iy$, where $\b_n=b_{n}+b_{n+1}+...+b_p$ and
$\e_n=a_na_{n+1}\cdot \cdot a_{p}$. Moreover,
\[
\lb{f01-4}
\p_{p-n}^+(z)={z^{p-n}\/\e_{p-n}}(1-z\b_{p-n+1}+O(z^{2})),\qqq
\p_{1}^+(z)=z{1-z\b_{2}+O(z^{2})\/\e_1}\qq as \qq z\to 0,
\]
\[
\lb{f01-5}
 \p_1^-(z)=z^{-1},\qqq \p_2^-(z)={1-zb_1\/a_1z^2}.
\]
\end{lemma}

\begin{lemma}
\lb{fw1} Let $q\in \ell_{1,2p}^2$.  Then
\[
\lb{fw1-1}
 w(z)=-(b_1-z^{-1})\p_1^+(z)-a_1\p_2^+(z),\qqq z\ne 0,
\]
\[
\lb{fw1-2}
 s(z)=(1-b_1z^{-1})\p_1^+(z)-z^{-1}a_1\p_2^+(z),
\]
\begin{multline}
\lb{fw1-3}
 w(z)={z^{2p-1}\/\e_1}\rt(-b_1c_p+
 {c_p(c_1+\b_{2}b_1)+b_1b_pa_p^2\/z}+{O(1)\/z^{2}}\rt)\qq
 as \qq z\to \iy\\
  w(z)={1-z\b_1+O(z^2)\/\e_1}\qqq  as \qq z\to0,
\end{multline}
\begin{multline}
\lb{fw1-4}
 s(z)={z^{2p-1}\/\e_1}\rt(c_p-{c_p\b_1+b_pa_p^2\/z}+{O(1)\/z^{2}}\rt)
\qqq   as \qq z\to \iy\\
  s(z)={1\/\e_1}\rt(-b_1+z(c_1+b_1\b_2)+O(z^2)\rt) \qqq  as \qq z\to 0.
\end{multline}
Furhermore, if $b_1=c_p=0$, then
\[
\lb{fw1-5}
 w(z)=-{z^{2p-3}\/\e_1}(c_1b_p+O(1/z))\qq
 as \qq z\to \iy.
\]
\end{lemma}
\no{\bf Proof.} Using \er{6} and Lemma \ref{f01}, we
obtain \er{fw1-1}, \er{fw1-2}. Then asymptotics from Lemma \ref{f01}
imply \er{fw1-3}, \er{fw1-4}. If $b_1=c_p=0$, then \er{f01-2}, \er{f01-3} imply
$$
w=-a_1\p_2^++{\p_1^+\/z}={a_1^2z^{2p-3}\/\e_1}(b_p+O(z^{-1}))-{z^{2p-3}\/\e_1}(b_p+O(z^{-1}))
=-{z^{2p-3}\/\e_1}(c_1b_p+O(z^{-1})).
$$
\BBox

\no{\bf Proof Proposition \ref{sw3}.} Identities \er{sw3-1} follow from \er{5}, \er{7}.
Recall the following identity
\[
A'(z_j)=-{1\/z_j}\sum _{n\in\Z}\p_n^+(z_j)\p_n^-(z_j),\qq j\in\Z_N,
\]
see (10.34) in \cite{Te}.
Then using \er{sw3-1} we obtain \er{sw3-2}. Similar arguments give
\er{sw3-3}.

Using $w(0)>0$ and \er{sw3-2}, we obtain  $w'(z_{1})<0$ and
$w'(z_{2})>0,...$. Moreover, due to \er{sw3-1}, \er{sw3-2} we have
\er{sw3-4}. Identity \er{fw1-1}, \er{fw1-2} give \er{sw3-5}.

Substituting $w=\sum_1^{2p}\check w_n z^{n-1}$ and $s=\sum_1^{2p}\check s_nz^{n-1}$ into \er{7} we obtain \er{trp}.
\BBox


We need some results  about the inverse problems from \cite{Te} for $q\in \ell_1^1=\{h=(h_n)_{n\in\Z}:
\sum (1+|n|)|h_n|<\iy\}$.
Define the Marchenko operator $\gF_n: \ell^2(\Z_+)\to \ell^2(\Z_+)$ by
\[
\lb{GLM1}
(\gF_nf)_k=\sum_{m\ge 0}F(2n+m+k)f_m, \qqq f=(f_n)_0^\iy\in
\ell^2(\Z_+),\qq \Z_+=\Z\cap [0,\iy),
\]
$n,k\ge 0$ where
\[
\lb{GLM2}
F(n)=\wh R_n^++\sum_{j=1}^N{z_j^n\/m_j},\qqq \wh R_n^\pm={1\/2\pi
i}\int_{|z|=1}R_\pm(z)z^{n-1}dz,\qqq n\in \Z,
\]
and  $\wh R_n^\pm$ are the Fourier coefficients of $R_\pm$.
The Marchenko equation is given by
\[
\lb{GLM3}
(I+\gF_n)\cK_n(\cdot)=\e_n^2e_0,\qq where  \qq \cK_n=(\cK_n(m))_0^\iy,\qq
e_n=(\d_{n,k})_0^\iy\in \ell^2(\Z_+),
\]
i.e.,
\[
\lb{GLM4}
\cK_n(k)+F(2n+k)+\sum_{m\ge 1}F(2n+k+m)\cK_n(m)=\e_n^2\d_{0,k}.
\]
For each $n\in\Z$ these equations have unique
"decreasing" solutions $(\cK_n(k))_{k=0}^\iy$. The sequences $a_n,b_n, n\in\Z$
have the forms
\[
\lb{recab} a_n^2={\P_n^0\/\P_{n+1}^0}, \ \ \
b_n={\P_{n}^1\/\P_{n}^0}- {\P_{n-1}^1\/\P_{n-1}^0} ,\qq
\P_{n}^k=\langle e_k,(I+\gF_{n})e_0\rangle, \qq k=0,1, \qqq n\in \Z,
\]
where $\langle \cdot, \cdot\rangle$ is the scalar product in
$\ell^2(\Z_+)$.
Recall the inverse spectral theorem from \cite{Te}.

\no \begin{theorem} \label{Tist}
The Faddeev mapping
\[
\ell_1^1(\Z)\os \ell_1^1(\Z)\to \gS=\{ R_{+}, \ (z_j,m_j^+, j\in
\Z_N)\in (-1,1)^N\ts \R_+^N, z_j\ne 0,\ N\ge 0\}
\]
given by $(a_n-1, b_n)_{n\in \Z}\to \{R_+(z),\ z\in \S^1,\  (z_j,
m_j^+, j\in \Z_N)\}$ is one-to-one and onto, where  the function
$R_{+}$, the eigenvalues $z_j$ and the norming constants $m_j^+,
j\in \Z_N$ satisfy the conditions

\no 1) $R_{+}(z)=\ol R_{+}(\ol z)=-{s(\ol z)\/zw(z)}, z\in  \S_0=\S^1\sm \{\pm 1\}$,
and the function $R_{+}(z)$ is continuous in $z \in\S_0$ and
satisfies
\[
\lb{4.2}
 C|1-z^2|^2+|R_{+}(z)|^2\le 1,\ \ \  all \ z \in\S_0, \ \ \ {\rm
for \ some} \ \ C>0.
\]
\no 2) The eigenvalues $z_j, j\in \Z_N$ are distinct and
$m_j^-m_j^+=w'(z_j)^2$.

\no 3) The sequences $\wh R^{\pm}=(\wh R_n^{\pm})_1^\iy$ defined in
\er{GLM2} with $R_{-}(z)=-R_+(\ol z){A(\ol z)\/A(z)}$ satisfy
\[
\sum_{n\ge 1}n|\wh R_n^{\pm}-\wh R_{n+2}^{\pm}|<\iy.
\]
\end{theorem}

We are ready to prove the first result.

\no {\bf Proof of Theorem \ref{T1}.} We consider the case $q\in
\gX_0^0$ and $m=2p-1$, the proof of other cases is similar.
 If $q\in  \gX_0^0$, then Lemma \ref{sw3} gives that  $(s,E_N)\in
\cS_0^0$, which yields a mapping $q\to (s,E_N)$ from $\gX_0^0$ into $\cS_0^0$.

 We will show uniqueness. Let $q\in  \gX_0^0$.
 Then Lemma \ref{sw3} gives that  $(s,E_N)\in
\cS_0^0$, and Theorem \ref{Ta1} gives a unique $w$. Moreover, Proposition \ref{sw3} yields the norming constants $m_j^{\pm}, j\in \Z_N$. These data determine the compactly supported potential uniquely by Theorem \ref{Tist}. Then we deduce that the mapping $q\to (s,E_N)$ is an injection.

We will show surjection of the mapping $q\to(s,E_N)$.
If $(s,E_N)\in \cS_0^0$, then Theorem \ref{Ta1}  gives unique
$w\in \cW_0^0$ and we have $R_+=-{s(z^{-1})\/zw(z)}$.
If $n\ge 2p+1$, then we have
$$
\wh R_n^+={1\/2\pi i}\int_{|z|=1}\!\!\!\!\!\!
R_+(z)z^{n-1}dz=-{1\/2\pi i}\int_{|z|=1}{z^{n-2}\wt
s(z)\/w(z)}dz=-\sum_{j=1}^N \Res {z^{n-2}\wt
s(z)\/w(z)}\rt|_{z_j}=-\sum_{j=1}^N{z_j^n\/m_j},
$$
since $n-2=(2p-1)+(n-2p-1)$ and the function $z^{2p-1}\wt s$ is a
polynomial. Thus
\[
F(n)=0\qq  \qq if \qqq n\ge 2p+1,
\]
and using \er{recab} we obtain
\[
\P_n^0=1+F(2n), \qq and\qq  a_n^2={1+F(2n)\/1+F(2n+2)}
\]
thus $a_n^2=1$ if $n>p$. Similar arguments yield $b_n=0$ if
$n>p$. Then we deduce that
$$
a_n=1,\qq b_n=0  \qq if \qqq n\ge p+1.
$$

Moreover, similar arguments yield $b_n=0, a_n=$ if $n\le 0$.
Then the asymptotics from Lemma \ref{fw1} give that $q\in \gX_0^0$, which
yields surjection.  \BBox

{\bf Proof of Theorem \ref{T2}.} i)  We consider the case $q\in \gX_0^0$, the proof of other cases is similar.
Let $q\in \gX_0^0$. Then Proposition \ref{sw3} yield
$w(\cdot,q)\in \cW_0^0$ and we have the mapping $q\to (w,\s)$, where the sequence $\s\in \X_w$ is given by
\er{setX0}, since $(s, E_N)\in \cS_0^0$.

By Theorem \ref{Ta2}, for each $(w,\s)\in (w,\X_w)$ there exists a unique $(s, E_N)\in \cS_0^0$,
then due to Theorem \ref{T1} the mapping  $q\to (w,\s)$ is an injection.

We will prove that $\gJ_R$ is a surjection. Let $w\in
\cW_0^0$ have zeros $E_N=(z_j)_{j\in\Z_N}$ from $(-1,1), N\geq 0$ and let a sequence $\s=(\s_n)_0^m \in\X_w$ be defined by \er{setX0}.
By Theorem \ref{Ta1}, there exists a unique  $(s,E_N)\in \cS_0^0$
 such that \er{S2} with $E_N$
hold true. Then $(s,E_N)\in \cS_0^0$ and, by Theorem \ref{T1}, there
exists a unique potential $q\in \gX_0^0$ with the scattering
function $s(\cdot)$. Thus the mapping $\gJ_{res}: \gX_\n^\t\to((w,\s)\in
\cW_\n^\t \ts\X_w)$ given by $q\to (w,(\s_n)_0^m)$  is one-to-one and onto, where
 $m=2p-1$.

ii) Furthermore, using similar arguments we deduce that
if we fix $q\in  \gX_\n^\t$, then the mapping $\P: r\to \s(s(r))$
is a bijection between the iso-resonance set of potentials $\Iso(q)$ and the set
of sequences $\X_w, w=w(\cdot,q)$.
\BBox

\section {Properties of polynomials $s,w$}
\setcounter{equation}{0}

In this section we will get the needed results about polynomials  $s,w$. In
order to solve the inverse problems $(a,b)\to $ (spectral data) we
need the following results.

\begin{lemma} \lb{sw5} Let $g=C_g\prod_1^m(z-\r_n)$ and $f=C_f\prod_1^m(z-\z_n)
$ satisfy
\[
\lb{A1}
 g(z)g(z^{-1})+\e^2=f(z)f(z^{-1}),\qq \e=z-z^{-1},\ all \ z\ne0,
\]
 for some $m\ge3$ and $g(0)\ne 0, f(0)\ne 0$. Then
\[
\lb{sw5-1} \prod_1^m(\l-\m_n)={\l^2-4\/C}+\prod_1^m(\l-\l_n),\qq
\m_n=\z_n+{1\/\z_n},\ \l_n=\r_n+{1\/\r_n},
\]
\[
\lb{sw5-2}
C=C_s^2C_\z=C_w^2C_\r,\qqq  \qq  {C_\r\/C_\z}>0,\qq C_\r=\prod_1^m(-\r_n),
\qq C_\z=\prod_1^m(-\z_n),
\]
\[
\lb{sw5-3}
\sum_1^{m}\m_n=\sum_1^{m}\l_n,\qq ...,\qq
\prod_1^m(-\m_n)+{4\/C}=\prod_1^m(-\l_n).
\]
\end{lemma}
\no{\bf Proof.} Using $f=C_f\prod_1^m (z-\z_n)$ and
$g=C_g\prod_1^m (z-\r_n)$, we obtain
$$
C_s^2\prod_1^m (z-\z_n)(z^{-1}-\z_n)=\e^2+C_w^2\prod_1^m
(z-\r_n)(z^{-1}-\r_n).
$$
Thus the identities $-\z_n(\l-\m_n)=(z-\z_n)(z^{-1}-\z_n)$ and
$-\r_n(\l-\l_n)=(z-\r_n)(z^{-1}-\r_n)$,  give
$$
C_s^2C_\z\prod_1^m (\l-\m_n)=\l^2-4+C_w^2C_\r\prod_1^m (\l-\l_n),\qq
$$
and $C_s^2C_\z=C_w^2C_\r=C$, which yields
${C_\z\/C_\r}={C_w^2\/C_s^2}>0$. \BBox

Assume that  we know only the polynomial $f$ in the equation \er{A1}
 and we have to determine $g$. In order to do this we have to solve
the equation \er{A1} in some class of polynomials. The following
Theorem will be used to determine the polynomial $w$ if we know $s$.

 \begin{theorem}
 \label{Ta1}
 i) Let functions $f,g$ be analytic in $\C\sm\{0\}$ and satisfy \er{A1}
 and be real on $\R\sm\{0\}$.
 Then $f^2(\pm 1)=g^2(\pm 1)$. Moreover, if
 $f(\pm 1)=0$ then $(f')^2(\pm 1)=8+(g')^2(\pm
 1)\ge 8$.

\no  ii) Let $(f,K_f^N)\in \cS_\n^\t, $ for some $(\n,\t)\in
\{0,1\}, N\ge 0, m\ge 3$. Then there exists a unique polynomial
$g\in \cW_\n^\t$ satisfying \er{A1}.
\end{theorem}
\no {\bf  Remark}. It is possible that the function $f(z)f(z^{-1})-\e^2(z)$ has more
zeros on $(-1,1)\sm \{0\}$, and it is important that we make a
special choice which, however, has to satisfy condition \er{defW2}.

\no {\bf Proof.} The statement i) is very simple and differentiating
\er{A1}  we obtain $(f')^2(\pm 1)=8+(f')^2(\pm 1)$ at $z=\pm 1$.

ii) Recall that $f=z^\n C\prod_{1}^{m} (z-\z_n)$, where
$$
C\in\R\sm\{0\}, \qq 0<|\z_{1}|\le |\z_{2}|\le ..\le |\z_{m}|, \qq
m=2p-1-\t-\n,
$$
and $K_f^N=(z_j)_{j\in \Z_N}$ is some sequence of zeros of
the function $f(z)f(z^{-1})-\e^2(z)$ such that
\[
\lb{a2} -1\!\!<z_{n_-}\!\!<...<\!\!z_{-1}\!\!<\!\!0
\!\!<\!\!z_{1}\!\!<...<\!\!z_{n_+}\!\!<\!\!1,\qq
 (-1)^{n_\pm}f(\pm1)\ge 0,\ \ z_j(-1)^{j}f(z_j^{\pm 1})>0,\ \
\]
$ j\in \Z_N$. Then
$$
f(z)f(z^{-1})=C^2\prod_{1}^{m}
(z-\z_n)(z^{-1}-\z_n)={C_0\/z^{m}}\prod_{1}^{m}
(z-\z_n)(z-\z_n^{-1}),\
$$
where $C_0=C^2C_\z, C_\z=\prod_{1}^{m}(-\z_n)$. Then
$G(z)=f(z)f(z^{-1})-\e^2(z)$ satisfies
\[
G(z)=\ca C_0z^{m}(1+O(z^{-1})) & as \qq z\to \iy \\
         C_0z^{-m}(1+O(z)) & as \qq z\to0 \ac,
\]
and thus
\[
\lb{eiG} G(z)=-\e^2(z)+{C_0\/z^{m}}\prod_{1}^{m}
(z-\z_n)(z-\z_n^{-1})
={C_0\/z^{m}}\prod_1^{m} (z-\r_n)(z-\r_n^{-1}),
\]
where $\r_n\ne 0, \r_n^{-1}$ are the zeros of $G$ counted with multiplicity and
satisfying
\begin{multline}
\lb{Az} 0<|\r_{1}|\le |\r_{2}|\le ..\le |\r_{N}|< 1\le
|\r_{N+1}|\le..\le |\r_{m}|, \\
the\ set \ \{\r_{1},\r_{2}, ...,\r_{N}\}=\{z_j, j\in \Z_N\}\ss
(-1,1)\sm\{0\},
\end{multline}
and Conditions i)-iv) in Definition W,   since $f$ satisfies
Definition S and $G(z)=G(1/z)$ for all $z\ne 0$. Moreover, using
$g_0(z)=\prod_1^{m}(z-\r_n)$ we have
\[
G(z)=Cg_0(z)g_0(z^{-1}),\qq where \qqq C={C_0\/g_0(0)},\qq
\]
Note that Lemma \ref{sw5}  gives $C>0$, then $g=C_*g_0$ and
$g(0)>0$, where $C_*$ satisfies $C_*g_0(0)>0, C_*^2=C>0$. By the
construction of $g$, this function is unique.
 $\BBox$

In order to solve the inverse problems $\gX_\n^\t\to  \cW_\n^\t$ we
need the following results. Assume that we have $(w,\s)$, where  the function $w\in
\cW_\n^\t$ and the sequence $\s\in\X_w$, then we have to determine $s$
uniquely. In order to do this we have to solve the equation \er{A1}
 in  class of polynomials $s\in \cS_\n^\t$. The sequence
$\s$ will give uniqueness. The following Theorem will be used to
determine the function $s$ if we know $(w,\s)$.

\begin{theorem}
 \label{Ta2}
Let $g\in \cW_\n^\t, $ for some $\n,\t\in \{0,1\}, m\ge 3$. Then for
each $\s\in \X_g$ defined in \er{setX0}  there exists a unique
$(f,K_f^N)\in \cS_\n^\t$ satisfying \er{A1}.

\end{theorem}

\no {\bf Proof.} Recall that $g=C_g\prod_1^m (z-\r_n), m=2p-1-\t-\n$
for some $C_g\in \R\sm\{0\}$ and $\r_n \in \C\sm\{0\}$ such that:

\no i) $g$ is real on $\R$ and if $+1$ and /or $-1$ are
zeros, they are simple, and $g(0)>0$,

\no ii) $|g(z)|\ge |\e(z)|$ for  any $|z|=1$, where $\e=z-{1\/z}$,

\no iii)  $g$ has only simple zeros $z_{n_-},..,z_{-1},
z_{1},..,z_{n_+}$ in $\dD_1$ for some $\pm n_\pm \ge 0$ such that
$-1<z_{n_-}<...<z_{-1}<0<z_{1}<...<z_{n_+}<1$.

Then we obtain
$$
g(z)g(z^{-1})=C_g^2\prod_1^{m}
(z-\r_n)(z^{-1}-\r_n)={C_0\/z^{m}}\prod_1^{m}
(z-\r_n)(z-\r_n^{-1}),\
$$
where $C_0=C_g^2\prod_1^{m}(-\r_n)$. Then $F(z)=\e^2(z)+g(z)g(z^{-1})$
satisfies
\[
F(z)=\ca C_0z^{m}(1+O(z^{-1})) & as \qq z\to \iy \\
         C_0z^{-m}(1+O(z)) & as \qq z\to0 \ac,
\]
and thus
\[
F(z)={C_0\/z^{m}}\prod_1^{2m} (z-t_n)={C_0\/z^{m}}\prod_1^{m}
(z-t_n)(z-t_n^{-1}),
\]
where $t_n\ne 0$ are the zeros of $F$ counted with multiplicity and
satisfying
$$
0<|t_{1}|\le |t_{2}|\le ..\le |t_{m}|, \qq t_{n+m}=1/t_n,\ all \
n=1,2,..,m, \qq \arg t_n\in [0,2\pi),
$$
where if $|t_n|=|t_k|,\ \arg t_n\le \arg t_k \Rightarrow n\le k$ and if $|t_n|=1$, then $\Im t_n\ge 0$, since $F(z)=F(1/z)$ for all $z\ne 0$.
 Moreover, we have
\[
F(z)={C_0\/z^{m}}\prod_1^{m}
(z-\z_n)(z-\z_n^{-1})={C_0\/f_0(0)}f_0(z)f_0(z^{-1}), \qq
  f_0(z)=\prod_1^{m} (z-\z_n),
\]
where  $\z_n=t_n^{\s_n},\qq (\s_n)_1^m\in \X_g$. Note that Lemma
\ref{sw5} gives $E={C_0\/f_0(0)}>0$, then $f=C_*f_0$ and $g(0)>0$,
where $C_*=\sqrt E$. We need to chose the sign of $C_*$. By the
construction of $g$, this function is unique.  $\BBox$

{\bf Example p=1}. In this case we have
\[
w={1\/a_1}(1-b_1z-a_1^2z^2),\qqq s={1\/a_1}(-b_1+c_1z).
\]
1) If we know $s=s_0+2s_1z$, then
\[
a_1=-s_1+\sqrt{s_1^2+1},\qq b_1=-{s_0\/s_1}.
\]
2) If we know $w$, then
\[
-a_1w=z^2+{b_1\/a_1^2}z-{1\/a_1^2},\qq
\r_1={-b_1-\sqrt{b_1^2+4a_1^2}\/4a_1^2}<0, \qq
\r_2={-b_1+\sqrt{b_1^2+4a_1^2}\/4a_1^2}>0,
\]
where $\r_1,\r_2\in \R$. Then
$$
\r_1\r_2=-{1\/a_1^2},\qq b_1={\r_1+\r_2\/\r_1\r_2}.
$$
{\bf Example p=2}. In this specific case we have
\[
\p_2^+={z^2\/a_2},\qq \p_1^+={c_2z^3-b_2z^2+z\/a_1a_2},
\]
and
\[
w={1\/a_1a_2}\rt(-z^3c_2b_1+z^2(c_2-a_1^2+b_1b_2)-z(b_1+b_2)+1
\rt),
\]
\[
s={1\/a_1a_2}\rt(z^3c_2-z^2(b_2+b_1c_2)+z(c_1+b_1b_2)-b_1  \rt).
\]
We obtain 4 cases:

\no 1) The case $a_2\ne 1, b_1\ne 0$ have been considered since
$m=3$.

\no 2) Let $a_2=1, b_2\ne 1, b_1=0, a_1\ne 1$. Then it is similar to
the case $p=1$ and we have:
\[
w={-a_1^2z^2-zb_2+1\/a_1},\qqq s=z{-zb_2+c_1\/a_1}.
\]
3) If $a_2=1, b_2\ne 0, b_1\ne 0$, then
\[
w={z^2(b_1b_2-a_1^2)-z(b_1+b_2)+1\/a_1},\qqq
s={-z^2b_2+z(c_1+b_1b_2)-b_1\/a_1}.
\]
4) If $a_2\ne 1, b_1=0, a_1\ne 1$, then
\[
w={z^2(c_2-a_1^2)-zb_2+1\/a_1a_2},\qqq s=z{z^2c_2-zb_2+c_1\/a_1a_2}.
\]


\no {\bf Acknowledgments.}  The research was partially supported  by EPSRC grant EP/D054621.


\begin{thebibliography} {9999}\setlength{\itemsep}{-\parskip}
\footnotesize


\bibitem[BNW] {BNW} B. M. Brown; S. Naboko; R. Weikard,
The inverse resonance problem for Jacobi operators Bull. London
Math. Soc. 37(2005), 727–37.


\bibitem[C1] {C1} K. M. Case, The discrete inverse scattering problem in one
dimension, J. Math. Phys. 15(1974), 143–146.

\bibitem[C2] {C2} K. M. Case, On discrete inverse scattering problems. II, J. Math.
Phys. 14(1973), 916–920.



\bibitem[CC] {CC}  K. M. Case; S. C. Chiu The discrete version of the Marchenko
equations in the inverse scattering problem, J. Math. Phys. 14(1973),
1643–1647.

\bibitem[CK] {CK} K. M. Case and M. Kac, A discrete version of the inverse
scattering problem, J. Math. Phys. 14(1973), 594–603.

\bibitem[DS1] {DS1} D. Damanik; B. Simon, Jost functions and Jost solutions for Jacobi matrices, I. A necessary and sufficient condition for Szeg\"o asymptotics, Invent. Math. 165 (2006), no. 1, 1-50.

\bibitem[DS2] {DS2} D. Damanik; B. Simon,
Jost functions and Jost solutions for Jacobi matrices. II. Decay and analyticity.
Int. Math. Res. Not. 2006, Art. ID 19396, 32 pp.


\bibitem[DT] {DT} Deift, P., Trubowitz, E. Inverse scattering on
the line, Commun. Pure and Applied Math., 32(1979), 121-251.

\bibitem[Fa] {Fa} Faddeev L., Properties of the S-matrix of the
one-dimensional Schr\"odinger equation, Trudy Mat. Inst. Steklov 73(1964), 314-333,
 English translation in AMST, 265, 139-166.


\bibitem[F] {F}  R. Froese, Asymptotic distribution of resonances in
one dimension, Journal of Differential Equations 137 (1997), no. 2, 251–272.


\bibitem[G1] {G1} G. S. Guseinov, The determination of an infinite Jacobi matrix
from the scattering data, Soviet Math. Dokl., 17(1976), 596–600.


\bibitem[G2] {G2}G. S. Guseinov, The inverse problem of scattering theory for a
second-order difference equation on the whole axis, Soviet Math.
Dokl., 17(1976), 1684–1688.





\bibitem[K1] {K1} E. Korotyaev,  Inverse resonance scattering on the real
line. Inverse Problems 21 (2005), no. 1, 325--341.

\bibitem[K2] {K2} E. Korotyaev, Inverse resonance scattering on the half line.
 Asymptot. Anal. 37 (2004), no. 3-4, 215--226.


\bibitem[K3] {K3} E. Korotyaev, Stability for inverse resonance problem.
Int. Math. Res. Not. 2004, no. 73, 3927--3936.

\bibitem[MBO] {MBO}
F. G. Maksudov, Eh. M. Bajramov, R. U. Orudzheva, The inverse
scattering problem for an infinite Jacobi matrix with operator
elements, Russ. Acad. Sci., Dokl., Math. 45(1992), No.2, 366–370.

\no  \bibitem[M] {M} Marchenko V. Sturm-Liouville operator and applications.
Basel: Birkh\"auser 1986.



\bibitem[MW] {MW} M. Marletta; R. Weikard, Stability for the inverse resonance problem
for a Jacobi operator with complex potential, Inverse Problems 23
(2007), 1677–1688.

\no  \bibitem[Me] {Me}  Melin A. Operator methods for inverse scattering on the real line.
Comm. P.D.E. 10(1985), 677-786.


\bibitem[NMPZ] {NMPZ}  S. Novikov; S. Manakov; L. Pitaevski; V. Zakharov,
Theory of solitons. The inverse scattering method. Consultants
Bureau [Plenum], New York, 1984.


\bibitem[S] {S} B. Simon, Resonances in one dimension and Fredholm
determinants Journal of Functional Analysis 178 (2000), no. 2, 396–420.

\bibitem[Te] {Te} G. Teschl, Jacobi operators and completely integrable
nonlinear lattices. Mathematical Surveys and Monographs, 72. AMS,
Providence, RI, 2000.

\bibitem[T1] {T1} M. Toda, Theory of Nonlinear Lattices, 2nd enl. ed., Springer,
Berlin, 1989.

\bibitem[T2] {T2} M. Toda, Theory of Nonlinear Waves and Solitons, Kluwer,
Dordrecht, 1989.

 \bibitem[Z] {Z} M. Zworski, Distribution of poles for scattering on
 the real line, Journal of Functional Analysis 73 (1987), no. 2, 277–296.

 \bibitem[Z1] {Z1} M. Zworski, A remark on isopolar potentials, SIAM Journal on Mathematical Analysis 32 (2001), no. 6, 1324–1326.

\end{thebibliography}
\end{document}